\documentclass[12pt]{article}
\usepackage{amsmath,amsfonts}
\numberwithin{equation}{section}
\newtheorem{theorem}{Theorem}[section]
\newtheorem{lemma}[theorem]{Lemma}

\newcommand{\qed}{\nolinebreak\hfill\vbox{\hrule\hbox{\vrule\kern3pt\vbox{\kern6pt}\kern3pt\vrule}\hrule}}
\newenvironment{pf}{{\noindent\bf Proof.}}{\qed\newline}
\newcommand{\plap}{\bigtriangleup_p}

\begin{document}

\title{The One-Phase Bifurcation For The $p$-Laplacian}
\author{Alaa Akram Haj Ali \& Peiyong Wang\footnote{Peiyong Wang is partially supported by a Simons Collaboration Grant.}\\
\footnotesize Department of Mathematics\\
\footnotesize Wayne State University\\
\footnotesize Detroit, MI 48202\\
\normalsize}
\date{}
\maketitle
\begin{abstract}
A bifurcation about the uniqueness of a solution of a singularly perturbed free boundary problem of phase transition associated with the $p$-Laplacian, subject to given boundary condition is proved in this paper. We show this phenomenon by proving the existence of a third solution through the Mountain Pass Lemma when the boundary data decreases below a threshold. In the second part, we prove the convergence of an evolution to stable solutions, and show the Mountain Pass solution is unstable in this sense.
\end{abstract}



\section{Introduction}\label{introduction}
In this paper, one considers the phase transition problem of minimizing the $p$-functional
\begin{equation}\label{p-functional}
J_{p,\varepsilon}(u) = \int_{\Omega}\frac{1}{p}|\nabla u(x)|^p + Q(x)\Gamma_{\varepsilon}(u(x))\,dx\ \ \ (1<p<\infty)
\end{equation}
which is a singular perturbation of the one-phase problem of minimizing the functional associated with the $p$-Laplacian
\begin{equation}\label{p-functional_original}
J_p(u) = \int_{\Omega}\frac{1}{p}|\nabla u(x)|^p + Q(x)\chi_{\{u(x)>0\}}\,dx,
\end{equation}
where $\Gamma_{\varepsilon}(s) = \Gamma(\frac{s}{
\varepsilon})$ for $\varepsilon > 0$ and for a $C^{\infty}$ function $\Gamma$ defined by
$$\Gamma(s) = \left\{\begin{array}{ll}
0 &\ \text{\ if\ }s\leq 0\\
1 &\ \text{\ if\ }s\geq 1,
\end{array}\right.$$
and $0\leq\Gamma(s)\leq 1$ for $0<s<1$, and $Q\in W^{2,2}(\Omega)$ is a positive continuous function on $\Omega$ such that $\inf_{\Omega}Q(x) > 0$. Let $\beta_{\varepsilon}(s) = \Gamma'_{\varepsilon}(s) = \frac{1}{
\varepsilon}\beta(\frac{s}{\varepsilon})$ with $\beta = \Gamma'$. The domain $\Omega$ is always assumed to be smooth in this paper for convenience. As in the following we will fix the value of $\varepsilon$ unless we specifically examine the influence of the value of $\varepsilon$ on the critical boundary data and will not use the notation $J_p$ for a different purpose, we are going to abuse the notation by using $J_p$ for the functional $J_{p,\varepsilon}$ from now on.

The Euler-Lagrange equation of (\ref{p-functional}) is
\begin{equation}\label{eulereq}
-\plap u + Q(x)\beta_{\varepsilon}(u) = 0\ \ x\in\Omega
\end{equation}
One imposes the boundary condition
\begin{equation}\label{bdrycondition}
u(x) = \sigma(x),\ \ x\in\partial\Omega
\end{equation}
on $u$, for $\sigma\in C(\partial\Omega)$ with $\min_{\partial\Omega}\sigma > 0$, to form a boundary value problem.

In this paper, we take on the task of establishing in the general case when $p\ne 2$ the results proved in \cite{CW} for the Laplacian when $p=2$. The main difficulty in this generalization lies in the lack of sufficient regularity and the singular-degenerate nature of the $p$-Laplacian when $p\ne 2$. A well-known fact about $p$-harmonic functions is the optimal regularity generally possessed by them is $C^{1,\alpha}$ (e.\,g.\,\cite{E} and \cite{Le}). Thus we need to employ more techniques associated with the $p$-Laplacian, and in a case or two we have to make our conclusion slightly weaker. Nevertheless, we follow the overall scheme of approach used in \cite{CW}. In the second section, we prove the bifurcation phenomenon through the Mountain Pass Theorem. In the third section, we establish a parabolic comparison principle. In the last section, we show the convergence of an evolution to a stable steady state in accordance with respective initial data.

\section{A Third Solution}\label{thirdsolution}
We first prove if the boundary data is small enough, then the minimizer is nontrivial. More precisely, let $u_0$ be the trivial solution of (\ref{eulereq}) and (\ref{bdrycondition}), being $p$-harmonic in the weak sense, and $u_2$ be a minimizer of the $p$-functional (\ref{p-functional}), and
set $$\sigma_M = \max_{\partial\Omega}\sigma(x)\ \ \text{ and\ }\ \ \sigma_m = \min_{\partial\Omega}\sigma(x).$$ If $\sigma_M$ is small enough, then $u_0\neq u_2$.

In fact, we pick $u\in W^{1,p}(\Omega)$ so that
\begin{equation}
\left\{\begin{array}{ll}
u = 0 &\ \ \text{ in $\Omega_{\delta}$}\\
u = \sigma &\ \ \text{ on $\partial\Omega$,\ \ \ \  and }\\
-\plap u = 0 &\  \ \text{ in $\Omega\backslash\bar{\Omega}_{\delta}$,}\end{array}\right.
\end{equation}
where $\Omega_{\delta} = \{x\in\Omega\colon dist(x,\partial\Omega) > \delta\}$ and $\delta > 0$ is a small constant independent of $\varepsilon$ and $\sigma$ so that $\int_{\Omega_{\delta}}Q(x)\,dx$ has a positive lower bound which is also independent of $\varepsilon$ and $\sigma$. Using an approximating domain if necessary, we may assume $\Omega_{\delta}$ possesses a smooth boundary. Clearly,
\begin{equation*}
J_p(u_0) = \int_{\Omega}\frac{1}{p}|\nabla u_0|^p + Q(x)\,dx
\geq \int_{\Omega}Q(x)\,dx.
\end{equation*}
It is well-known that
\begin{equation*}
\int_{\Omega\backslash\Omega_{\delta}}|\nabla u|^p \leq C\sigma^{\,p}_M\delta^{1-p}\ \ \text{ for $C = C(n,p,\Omega)$},
\end{equation*}
so that
\begin{alignat*}{1}
&J_p(u) = \int_{\Omega\backslash \Omega_{\delta}}\frac{1}{p}|\nabla u|^p +
\int_{\Omega\backslash\Omega_{\delta}}Q(x)\,dx\\
&\leq C\sigma^{\,p}_M\delta^{1-p} +
\int_{\Omega\backslash\Omega_{\delta}}Q(x)\,dx.
\end{alignat*}
So, for all small $\varepsilon > 0$,
\begin{equation*}
J_p(u) - J_p(u_0) \leq C\sigma^{\,p}_M\delta^{1-p} - \int_{\Omega_{\delta}} Q(x)\,dx < 0
\end{equation*}
if $\sigma_M\leq \sigma_0$ for some small enough $\sigma_0 = \sigma_0(\delta, \Omega, Q)$.

Let $\mathfrak{B}$ denote the Banach space $W^{1,p}_0(\Omega)$ we will work with. For every $v\in\mathfrak{B}$, we write $u = v + u_0$ and adopt the norm $\|v\|_{\mathfrak{B}} = \left(\int_{\Omega}|\nabla v|^p\right)^{\frac{1}{p}} = \left(\int_{\Omega}|\nabla u - \nabla u_0|^p\right)^{\frac{1}{p}}$. We define the functional
\begin{equation}
I[v] = J_p(u) - J_p(u_0) = \int_{\Omega}\frac{1}{p}|\nabla u|^p - \int_{\{u < \varepsilon\}}Q(x)\left(1 - \Gamma_{\varepsilon}(u)\right) - \int_{\Omega}\frac{1}{p} |\nabla u_0|^p
\end{equation}
Set $v_2 = u_2 - u_0$. Clearly, $I[0] = 0$ and $I[v_2] \leq 0$ on account of the definition of $u_2$ as a minimizer of $J_p$. If $I[v_2] < 0$ which is the case if $\sigma_M$ is small, we will apply the Mountain Pass Lemma to prove there exists a critical point of the functional $I$ which is a weak solution of the problem (\ref{eulereq}) and (\ref{bdrycondition}).

The Fr\'{e}chet derivative of $I$ at $v\in \mathfrak{B}$ is given by
\begin{equation}
I'[v]\varphi = \int_{\Omega}|\nabla u|^{p-2}\nabla u\cdot\nabla\varphi + Q(x)\beta_{\varepsilon}(u)\varphi\ \ \ \ \varphi\in\mathfrak{B}
\end{equation}
which is obviously in the dual space $\mathfrak{B}^*$ of $\mathfrak{B}$ in light of the H\"{o}lder's inequality. Equivalently
\begin{equation}
I'[v] = -\plap (v+u_0) + Q(x)\beta_{\varepsilon}(v + u_0)\in\mathfrak{B}^*.
\end{equation}
We see that $I'$ is Lipschitz continuous on any bounded subset of $\mathfrak{B}$ with Lipschitz constant depending on $\varepsilon$, $p$, and $\sup Q$. In fact, for any $v$, $w$, and $\varphi\in \mathfrak{B}$,
\begin{alignat*}{1}
&\ \ \left|I'[v]\varphi - I'[w]\varphi\right| = |\int_{\Omega}|\nabla v + \nabla u_0|^{p-2}(\nabla v + \nabla u_0)\cdot\nabla\varphi + Q(x)\beta_{\varepsilon}(v+u_0)\\
&- |\nabla w + \nabla u_0|^{p-2}(\nabla w + \nabla u_0)\cdot\nabla\varphi - Q(x)\beta_{\varepsilon}(w+u_0)| \\
&\leq \left|\int_{\Omega}|\nabla v + \nabla u_0|^{p-2}(\nabla v + \nabla u_0)\cdot\nabla\varphi -
|\nabla w + \nabla u_0|^{p-2}(\nabla w + \nabla u_0)\cdot\nabla\varphi\right| \\
&+ \left|\int_{\Omega}Q(x)\beta_{\varepsilon}(v+u_0) - Q(x)\beta_{\varepsilon}(w+u_0)\right|
\end{alignat*}
Furthermore,
\begin{alignat*}{1}
&\ \ \ \ \left|\int_{\Omega}Q(x)\beta_{\varepsilon}(v+u_0) - Q(x)\beta_{\varepsilon}(w + u_0)\right|\\
&= \left|\int_{\Omega}Q(x)\int^1_0\beta'_{\varepsilon}((1-t)w + tv + u_0)\,dt\,(v(x) - w(x))\,dx \right|\\
&\leq \sup|\beta'_{\varepsilon}|\int_{\Omega}\left|Q(x)\left(v(x) - w(x)\right)\right|\,dx \\
&\leq \frac{C}{\varepsilon^2}\left(\int_{\Omega}Q^{p'}(x)\right)^{\frac{1}{p'}}\left(\int_{\Omega}|v(x) - w(x)|^p\,dx\right)^{\frac{1}{p}}
\end{alignat*}
and
\begin{alignat*}{1}
&\ \ \ \ \left|\int_{\Omega}|\nabla v + \nabla u_0|^{p-2}(\nabla v + \nabla u_0)\cdot \nabla \varphi - |\nabla w + \nabla u_0|^{p-2}(\nabla w + \nabla u_0)\cdot \nabla \varphi\right|\\
&\leq \left|\int_{\Omega}|\nabla v + \nabla u_0|^{p-2}(\nabla v - \nabla w)\cdot \nabla \varphi\right| \\
&\ \ \ \ + \left|\int_{\Omega}\left(|\nabla v + \nabla u_0|^{p-2} - |\nabla w + \nabla u_0|^{p-2}\right)(\nabla w + \nabla u_0)\cdot \nabla \varphi\right|.
\end{alignat*}
In addition,
\begin{alignat*}{1}
&\ \ \ \ \left|\int_{\Omega}|\nabla v + \nabla u_0|^{p-2}(\nabla v - \nabla w)\cdot \nabla \varphi\right|\\
&\leq \left(\int_{\Omega}|\nabla v + \nabla u_0|^p\right)^{\frac{p-2}{p}}\left(\int_{\Omega}|\nabla\varphi|^p\right)^{\frac{1}{p}}\left(\int_{\Omega}|\nabla v - \nabla w|^p\right)^{\frac{1}{p}},
\end{alignat*}
and
\begin{alignat*}{1}
&\ \ \ \ \left|\int_{\Omega}\left(|\nabla v + \nabla u_0|^{p-2} - |\nabla w + \nabla u_0|^{p-2}\right)\left(\nabla w + \nabla u_0\right)\cdot\nabla\varphi\right|\\
&\leq C(p)\int_{\Omega}\left(|\nabla v + \nabla u_0|^{p-3} + |\nabla w + \nabla u_0|^{p-3}\right)|\nabla v - \nabla w||\nabla w + \nabla u_0||\nabla\varphi|\\
&\leq C(p)\left(\|\nabla v\|_{L^p} + \|\nabla w\|_{L^p} + \|\nabla u_0\|_{L^p}\right)^{p-2}\|\nabla v - \nabla w\|_{L^p(\Omega)}\|\nabla\varphi\|_{L^p(\Omega)}.
\end{alignat*}
Therefore $I'$ is Lipschitz continuous on bounded subsets of $\mathfrak{B}$.

We note that $f\in\mathfrak{B}^*$ if and only if there exist $f^0$, $f^1$, $f^2$, ..., $f^n\in L^{p'}(\Omega)$, where $\frac{1}{p} + \frac{1}{p'} = 1$, such that
\begin{alignat}{1}
&<f,u>\ = \int_{\Omega}f^0u + \sum^n_{i=1}f^iu_{x_i} \ \ \text{ holds for all $u\in\mathfrak{B}$; and}\label{repre}\\
&\|f\|_{\mathfrak{B}^*} = \inf\left\{\left(\int_{\Omega}\sum^n_{i=0}|f^i|^{p'}\,dx\right)^{\frac{1}{p'}}\colon \text{(\ref{repre}) holds.}\right\}
\end{alignat}

Next we justify the Palais-Smale condition on the functional $I$. Suppose $\{v_k\}\subset\mathfrak{B}$ is a Palais-Smale sequence in the sense that
\begin{equation*}
\left|I[v_k]\right|\leq M\ \ \ \ \text{and\ \ }\ \ I'[v_k]\rightarrow 0\ \ \ \ \text{in $\mathfrak{B}^*$}
\end{equation*}
for some $M > 0$. Let $u_k = v_k + u_0\in W^{1,p}(\Omega)$, $k = 1, 2, 3, ...$.

That $Q(x)\beta_{\varepsilon}(v + u_0)\in W^{1,p}_0(\Omega)$ implies that the mapping $v\mapsto Q(x)\beta_{\varepsilon}(v + u_0)$ from $W^{1,p}_0(\Omega)$ to $\mathfrak{B}^*$ is compact due to the fact $W^{1,p}_0(\Omega)\subset\subset L^p(\Omega)\subset \mathfrak{B}^*$ following from the Rellich-Kondrachov Compactness Theorem. Then there exists $f\in L^p(\Omega)\subset\mathfrak{B}^*$ such that for a subsequence, still denoted by $\{v_k\}$, of $\{v_k\}$, it holds that
\begin{equation*}
Q(x)\beta_{\varepsilon}(u_k)\rightarrow -f\ \ \text{ in $L^p(\Omega)$.}
\end{equation*}
Recall that
\begin{equation*}
\left|I'[v_k]\varphi\right| = \sup_{\|\varphi\|_{\mathfrak{B}}\leq 1}\left|\int_{\Omega}|\nabla u_k|^{p-2}\nabla u_k\cdot\nabla\varphi + Q(x)\beta_{\varepsilon}(u_k) \varphi\right|\rightarrow 0.
\end{equation*}
As a consequence,
\begin{equation}\label{test1}
\sup_{\|\varphi\|_{\mathfrak{B}}\leq M}\left|\int_{\Omega}|\nabla u_k|^{p-2}\nabla u_k\cdot\nabla \varphi - f\varphi\right| \rightarrow 0\ \ \ \ \text{for any $M\geq 0$.}
\end{equation}
Obviously, that $\{I[v_k]\}$ is bounded implies that a subsequence of $\{v_k\}$, still denoted by $\{v_k\}$ by abusing the notation without confusion, converges weakly in $\mathfrak{B} = W^{1, p}_0(\Omega)$. In particular,
\begin{equation*}
\int_{\Omega}fv_k - fv_m\rightarrow 0\ \ \ \ \text{as $k$, $m\rightarrow\infty$.}
\end{equation*}
Then by setting $\varphi = v_k - v_m = u_k - u_m$ in (\ref{test1}), one gets
\begin{equation}\label{conv}
\left|\int_{\Omega}\left(|\nabla u_k|^{p-2}\nabla u_k - |\nabla u_m|^{p-2}\nabla u_m\right)\cdot \nabla (u_k - u_m)\right| \rightarrow 0\ \ \ \ \text{as $k$, $m\rightarrow\infty$,}
\end{equation}
since
\begin{equation*}
\|u_k - u_m\|^p_{\mathfrak{B}} = \|v_k - v_m\|^p_{\mathfrak{B}} \leq 2pM + 2J_p[u_0].
\end{equation*}
In particular, if $p = 2$, $\{v_k\}$ is a Cauchy sequence in $W^{1,2}_0(\Omega)$ and hence converges.
We will apply the following elementary inequalities associated with the $p$-Laplacian, \cite{L}, to the general case $p\neq 2$:
\begin{alignat}{1}
&<|b|^{p-2}b - |a|^{p-2}a,\,b - a> \geq (p-1)|b-a|^2(1 + |a|^2 + |b|^2)^{\frac{p-2}{2}},\ \ 1\leq p\leq 2;\label{ele1}\\
&\text{and}\ \ \ \ <|b|^{p-2}b - |a|^{p-2}a,\,b - a> \geq 2^{2-p}|b-a|^p,\ \ p\geq 2.\label{ele2}
\end{alignat}
We assume first $1 < p < 2$. Let $K = 2pM + 2J_p[u_0]$. Then the first elementary inequality (\ref{ele1}) implies
\begin{alignat*}{1}
&\ \ \ \ \ (p-1)\int_{\Omega}|\nabla u_k - \nabla u_m|^2\left(1+|\nabla u_k|^2 + |\nabla u_m|^2\right)^{\frac{p-2}{2}}\\ &\leq \int_{\Omega}\left(|\nabla u_k|^{p-2}\nabla u_k - |\nabla u_m|^{p-2}\nabla u_m\right)\cdot \nabla (u_k - u_m) \rightarrow 0
\end{alignat*}
Meanwhile H\"{o}lder's inequality implies
\begin{alignat*}{1}
&\ \ \ \ \ \int_{\Omega}|\nabla v_k - \nabla v_m|^p = \int_{\Omega}|\nabla u_k - \nabla u_m|^p \\
&\leq \left(\int_{\Omega}|\nabla u_k - \nabla u_m|^2\left(1 + |\nabla u_k|^2 + |\nabla u_m|^2\right)^{\frac{p-2}{2}}\right)^{\frac{p}{2}}
\left(\int_{\Omega}\left(1 + |\nabla u_k|^2 + |\nabla u_m|^2\right)^{\frac{p}{2}}\right)^{\frac{2-p}{2}} \\
&\leq C(p)\left(|\Omega| + K\right)^{\frac{2-p}{2}} \left(\int_{\Omega}|\nabla u_k - \nabla u_m|^2\left(1 + |\nabla u_k|^2 + |\nabla u_m|^2\right)^{\frac{p-2}{2}}\right)^{\frac{p}{2}}
\end{alignat*}
Therefore, $\{v_k\}$ is a Cauchy sequence in $\mathfrak{B}$ and hence converges.

Suppose $p > 2$. The second elementary inequality (\ref{ele2}) implies
\begin{alignat*}{1}
&\ \ \ \ \ \int_{\Omega}|\nabla v_k - \nabla v_m|^p = \int_{\Omega}|\nabla u_k - \nabla u_m|^p \\
&\leq 2^{p-2}\int_{\Omega}\left(|\nabla u_k|^{p-2}\nabla u_k - |\nabla u_m|^{ p-2}\nabla u_m\right)\cdot \left(\nabla u_k - \nabla u_m\right),
\end{alignat*}
which in turn implies $\{v_k\}$ is a Cauchy sequence in $\mathfrak{B}$ and hence converges, on account of (\ref{conv}). The Palais-Smale condition is verified for $1 < p < \infty$ for the functional $I$ on the Banach space $W^{1,p}_0(\Omega)$.

Before we continue the main proof, let us state an elementary result closely related to the $p$-Laplacian, which follows readily from the Fundamental Theorem of Calculus.
\begin{lemma}\label{p-inequalities}
For any $a$ and $b\in\mathbb{R}^n$, it holds
\begin{equation}\label{ele3}
|b|^p \geq |a|^p + p<|a|^{p-2}a, b-a> +\, C(p)|b - a|^p\ \ \ \ (p\geq 2)
\end{equation}
where $C(p) > 0$.

If $1 < p < 2$, then
\begin{equation}\label{ele4}
|b|^p \geq |a|^p + p<|a|^{p-2}a, b-a> +\, C(p)|b-a|^2\int^1_0\int^t_0\left|(1-s)a+sb\right|^{p-2}\,dsdt,
\end{equation}
where $C(p) = p(p-1)$.
\end{lemma}

We are now in a position to show there is a closed mountain ridge around the origin of the Banach space $\mathfrak{B}$ that separates $v_2$ from the origin with the energy $I$ as the elevation function, which is the content of the following lemma.
\begin{lemma}
For all small $\varepsilon > 0$ such that $C\varepsilon \leq \frac{1}{2}\sigma_m$ for a large universal constant $C$, there exist positive constants $\delta$ and $a$ independent of $\varepsilon$, such that, for every $v$ in $\mathfrak{B}$ with $\|v\|_{\mathfrak{B}} = \delta$, the inequality $I[v] \geq a$ holds.
\end{lemma}
\begin{pf}
It suffices to prove $I[v] \geq a > 0$ for every $v\in C^{\infty}_0(\Omega)$ with $\|v\|_{\mathfrak{B}} = \delta$ for $\delta$ small enough, as $I$ is continuous on $\mathfrak{B}$, and $C^{\infty}_0(\Omega)$ is dense in $\mathfrak{B}$.

Let $\Lambda = \{x\in\Omega\colon u(x)\leq\varepsilon\}$, where $u = v + u_0$. We claim that $\Lambda = \emptyset$ if $\delta$ is small enough. If not, one may pick $z\in\Lambda$. Let $\mathcal{AC}([a,b], S)$ be the set of absolutely continuous functions $\gamma\colon [a,b]\rightarrow S$, where $S\subseteq\mathbb{R}^n$. For each $\gamma\in\mathcal{AC}([a,b], S)$, we define its length to be $L(\gamma) = \int^b_a|\gamma'(t)|\,dt$. For $x_0\in\partial\Omega$, we define the distance from $x_0$ to $z$ to be
\begin{equation*}
d(x_0,z) = \inf\{L(\gamma): \gamma\in\mathcal{AC}([0,1],\bar{\Omega}),
\ \text{s.t.\ }\gamma(0) = x_0, \ \text{and\ }\gamma(1)=z\}
\end{equation*}
As shown in \cite{CW}, there is a minimizing path $\gamma_{x_0}$ for the distance $d(x_0, z)$.

Suppose the domain $\Omega$ is convex or star-like about $z$. For any $x_0\in\partial\Omega$, let $\gamma = \gamma_{x_0}$ be a minimizing path of $d(x_0, z)$.
Then it is clear that $\gamma$ is a straight line segment and $\gamma(t)\neq z$
for $t\in [0,1)$. Furthermore, for any two distinct points $x_1$ and $x_2\in\partial\Omega$, the
corresponding minimizing paths do not intersect in $\Omega\backslash\{z\}$. For this reason,
we can carry out the following computation. Clearly $v(x_0) = 0$
and $v(\gamma(1)) = \varepsilon - u_0(\gamma(1))\leq \varepsilon - \sigma_m < 0$. So the
Fundamental Theorem of Calculus
\begin{equation*}
v(\gamma(1)) - v(\gamma(0)) = \int^1_0\nabla v(\gamma(t))\cdot\gamma'(t)dt
\end{equation*}
implies
\begin{equation}\label{ineq-ftc}
\sigma_m - \varepsilon \leq \int^1_0|\nabla v(\gamma(t))||\gamma'(t)|dt.
\end{equation}
For each $x_0\in\partial\Omega$, let $e(x_0)$ be the unit vector in the direction of $x_0 - z$ and $\nu(x_0)$ the outer normal to $\partial\Omega$ at $x_0$. Then $\nu(x_0)\cdot e(x_0) > 0$ everywhere on $\partial\Omega$. Hence the above inequality (\ref{ineq-ftc}) implies
\begin{alignat*}{1}
&\ \ \ \ (\sigma_m - \varepsilon)\int_{\partial\Omega}\nu(x_0)\cdot e(x_0)\,dH^{n-1}(x_0)\\
&\leq \int_{\partial\Omega}\int^1_0|\nabla v(\gamma(t))| |\gamma'(t)|\,dt\,\nu(x_0)\cdot e(x_0)\,dH^{n-1}(x_0) \\
&\leq \int_{\partial\Omega}\left(\int^1_0|\gamma'(t)|\,dt\right)^{\frac{1}{p'}}\left(\int^1_0|\nabla v(
\gamma(t))|^p|\gamma'(t)|\,dt\right)^{\frac{1}{p}}\nu(x_0)\cdot e(x_0)\,dH^{n-1}(x_0),\\
&\hspace{3.5in}\text{ where $\frac{1}{p} + \frac{1}{p'} = 1$,} \\
&= \int_{\partial\Omega} L(\gamma_{x_0})^{\frac{1}{p'}}\left(\int^1_0|\nabla v(\gamma(t))|^p|\gamma'(t)|
\,dt\right)^{\frac{1}{p}}\nu(x_0)\cdot e(x_0)\,dH^{n-1}(x_0) \\
&\leq \left(\int_{\partial\Omega}L(\gamma_{x_0})\nu(x_0)\cdot e(x_0)\,dH^{n-1}\right)^{\frac{1}{p'}}\left(\int_{\partial\Omega}
\int^1_0|\nabla v(\gamma(t))|^p|\gamma'(t)|\nu \cdot e\,dt\,dH^{n-1}\right)^{\frac{1}{p}} \\
&= C|\Omega|^{\frac{1}{p'}}\left(\int_{\Omega}|\nabla v|^p\,dx\right)^{\frac{1}{p}}\\
&\leq C|\{u > \varepsilon\}|^{\frac{1}{p'}}\delta \leq C|\{u>0\}|^{\frac{1}{p'}}\delta,
\end{alignat*}
where the second and third inequalities are due to the application of the H\"{o}lder's
inequality, and the constant $C$ depends on $n$ and $p$. The second equality follows from the two representation formulas
\begin{equation*}
\left|\Omega\right| = C(n)\int_{\partial\Omega}L(\gamma_{x_0})\nu(x_0)\cdot e(x_0)\,dH^{n-1}(x_0)
\end{equation*}
and
\begin{equation*}
\int_{\Omega}\left|\nabla v(x)\right|^p\,dx = C(n)\int_{\partial\Omega}\int^1_0\left|\nabla v(\gamma_{x_0}(t))\right|^p\,\left|\gamma'_{x_0}(t)\right|\nu(x_0) \cdot e(x_0)\,dt\,dH^{n-1}(x_0).
\end{equation*}
If we take $\delta$ sufficiently small and independent of $\varepsilon$ in the preceding inequality
\begin{equation*}
(\sigma_m - \varepsilon)\int_{\partial\Omega}\nu(x_0)\cdot e(x_0)\,dH^{n-1}(x_0) \leq C|\{u>0\}|^{\frac{1}{p'}}\delta,
\end{equation*}
the measure $|\{u > 0\}|$ of the positive domain would be greater than that of $\Omega$, which is
impossible, provided that
\begin{equation}\label{direction-normal-ineq}
\int_{\partial\Omega}\nu(x_0)\cdot e(x_0)\,dH^{n-1}(x_0) \geq C,
\end{equation}
for a constant $C$ which depends on $n$, $p$ and $|\Omega|$, but not on $z$ or $v$.
Hence $\Lambda$ must be empty. So we need to justify the inequality (\ref{direction-normal-ineq}). To fulfil that condition, for $e = e(x_0)$, we set $l(e,z) = l(e) = L(\gamma_{x_0})$. Then \begin{equation*}
\int_{\partial\Omega}\nu(x_0)\cdot e(x_0)\,dH^{n-1}(x_0) = \int_{e\in \partial B}\left(l(e)\right)^{n-1}\,d\sigma(e),
\end{equation*}
where $B$ is the unit ball about $z$ and $d\sigma(e)$ is the surface area element on the unit sphere $\partial B$ which is invariant under rotation and reflection. Clearly,
\begin{equation*}
\left(\int_{\partial B}\left(l(e)\right)^{n-1}\,d\sigma(e)\right)^{\frac{2}{n-1}}\geq C(n)\int_{\partial B}l^2(e)\,d\sigma(e)
\end{equation*}
Consequently, in order to prove (\ref{direction-normal-ineq}), one needs only to prove
\begin{equation}\label{equiv-integ}
\int_{\partial B}l^2(e)\,d\sigma(e) \geq C(n, p, |\Omega|).
\end{equation}
Next, we show the integral on the left-hand-side of (\ref{equiv-integ}) is minimal if $\Omega$ is a ball while its measure is kept unchanged. In fact, this is almost obvious if one notices the following fact. Let $\pi$ be any hyperplane passing through $z$, and $x_1$ and $x_2$ be the points on $\partial\Omega$ which lie on a line perpendicular to $\pi$. Let $x^*_1$ and $x^*_2$ be the points on the boundary $\partial\Omega_{\pi}$, where $\Omega_{\pi}$ is the symmetrized image of $\Omega$ about the hyperplane $\pi$, which lie on the line $\overline{x_1x_2}$. Let $2a = |\overline{x_1x_2}| = |\overline{x^*_1x^*_2}|$ and $d$ be the distance from $z$ to the line $\overline{x_1x_2}$. Then for some $t$ in $-a \leq t \leq a$, it holds that
\begin{equation*}
L^2(\gamma_{x_1}) + L^2(\gamma_{x_2}) = \left(d^2+(a-t)^2\right) + \left(d^2+(a+t)^2\right) \geq 2(d^2 + a^2) = 2\left(L^*(\gamma_{x^*_1})\right)^2.
\end{equation*}
As a consequence, if $\Omega^*$ is the symmetrized ball with measure equal to that of $\Omega$, then
\begin{equation*}
\int_{\partial B}l^2(e)\,d\sigma(e) \geq \int_{\partial B}\left(l^*(e)\right)^2\,d\sigma(e) = C(n,|\Omega|),
\end{equation*}
where $l^*$ is the length from $z$ to a point on the boundary $\partial\Omega^*$ which is constant.
This finishes the proof of the fact that $\Lambda = \emptyset$.

In case the domain $\Omega$ is not convex, the minimizing paths of $d(x_1, z)$ and
$d(x_2, z)$ for distinct $x_1$, $x_2\in\partial\Omega$ may partially coincide. We form the set $\mathcal{DA}(\partial\Omega)$ of the
points $x_0$ on $\partial\Omega$ so that a minimizing path $\gamma$ of $d(x_0, z)$
satisfies $\gamma(t)\in\Omega\backslash\{z\}$ for $t\in (0,1)$. We call a point in $\mathcal{DA}(\partial\Omega)$ a \textbf{directly accessible} boundary point. Let $\Omega_1$ be the union of these minimizing paths for the directly accessible boundary points. It is not difficult to see that $|\Omega_1| > 0$ and hence $H^{n-1}(\mathcal{DA}(\partial\Omega)) > 0$. Then we may apply the above computation to the star-like domain $\Omega_1$ with minimal modification. We have
\begin{equation}
(\sigma_m - C\varepsilon)\int_{\partial\Omega}\nu(x_0)\cdot e(x_0)\,dH^{n-1}(x_0) \leq C|\Omega_1|^{\frac{1}{p'}}\delta \leq C|\Omega|^{\frac{1}{p'}}\delta.
\end{equation}
For small enough $\delta$, this raises a contradiction $|\Omega| > |\Omega|$.
So $\Lambda = \emptyset$.

Finally we prove that $\|v\|_{\mathfrak{B}} = \delta$ implies
\begin{equation}
I[v]= \int_{\Omega}\frac{1}{p}|\nabla v + \nabla u_0|^p - \frac{1}{p}|\nabla u_0|^p \geq a\ \ \text{for a certain $a > 0$.}
\end{equation}

If $p\geq 2$, then the elementary inequality (\ref{ele3}) implies that
\begin{alignat*}{1}
I[v] &= \int_{\Omega}\frac{1}{p}\left|\nabla v+ \nabla u_0\right|^p - \frac{1}{p}\left|\nabla u_0\right|^p \\
&\geq \int_{\Omega}<\left|\nabla u_0\right|^{p-2}\nabla u_0, \nabla v> + C(p)\left|\nabla v\right|^p \\
&= C(p)\int_{\Omega}\left|\nabla v\right|^p = C(p)\delta^p > 0,
\end{alignat*}
while if $1 < p < 2$, then the elementary inequality (\ref{ele4}) implies
\begin{alignat*}{1}
I[v] &\geq p(p-1)\int_{\Omega}\left|\nabla v\right|^2\int^1_0\int^t_0\frac{1}{\left|\nabla u_0 + s\nabla v\right|^{2-p}}\,dsdtdx \\
&\geq p(p-1)\int_{\Omega}\left|\nabla v\right|^2\int^1_0\int^t_0\frac{1}{\left(\left|\nabla u_0\right| + s\left|\nabla v\right|\right)^{2-p}}\,dsdtdx. 
\end{alignat*}
If $\int_{\Omega}|\nabla u_0|^p = 0$, then $I[v] = \frac{1}{p}\delta^p > 0$. So in the following, we assume $\int_{\Omega}|\nabla u_0|^p > 0$.

Let $S = S_{\lambda} = \{x\in\Omega\colon |\nabla v| > \lambda\delta\}$, where the constant $\lambda = \lambda(p,|\Omega|)$ is to be taken. Then
\begin{alignat*}{1}
\delta^p &= \int_{\Omega}|\nabla v|^p = \int_{\{|\nabla v|\leq \lambda\delta\}}|\nabla v|^p + \int_S|\nabla v|^p \\
&\leq (\lambda\delta)^p|\Omega| + \int_S|\nabla v|^p
\end{alignat*}
and hence
\begin{equation*}
\int_S|\nabla v|^p \geq \delta^p\left(1 - \lambda^p|\Omega|\right) \geq \frac{1}{2}\delta^p,\ \ \text{if $\lambda$ satisfies\ }\frac{1}{4} < \lambda^p|\Omega| \leq \frac{1}{2}.
\end{equation*}
Meanwhile, for $1 < p < 2$, it holds that
\begin{alignat*}{1}
I[v] &\geq C(p)\int_{S}\left|\nabla v\right|^2\int^1_0\int^t_0\frac{1}{\left(\left|\nabla u_0\right| + s\left|\nabla v\right|\right)^{2-p}}\,dsdtdx \\
&=C(p)\left(\int_{S\cap\{|\nabla u_0|\leq |\nabla v|\}}\left|\nabla v\right|^2\int^1_0\int^t_0\frac{1}{\left(|\nabla u_0| + s|\nabla v|\right)^{2-p}}\,dsdtdx \right. \\
&\ \ \ \ + \left.\int_{S\cap\{|\nabla u_0| > |\nabla v|\}}\left|\nabla v\right|^2\int^1_0\int^t_0\frac{1}{\left(\left|\nabla u_0\right| + s\left|\nabla v\right|\right)^{2-p}}\,dsdtdx\right).
\end{alignat*}
The first integral on the right satisfies
\begin{alignat*}{1}
&\ \ \ \ \int_{S\cap\{|\nabla u_0|\leq |\nabla v|\}}\left|\nabla v\right|^2\int^1_0\int^t_0\frac{1}{\left(|\nabla u_0| + s|\nabla v|\right)^{2-p}}\,dsdtdx \\
&\geq \int_{S\cap\{|\nabla u_0|\leq |\nabla v|\}}\left|\nabla v\right|^p\int^1_0\int^t_0\frac{1}{\left(1 + s\right)^{2-p}}\,dsdtdx \\
&= C(p)\int_{S\cap\{|\nabla u_0|\leq |\nabla v|\}}\left|\nabla v\right|^p\,dx,
\end{alignat*}
while the second integral on the right satisfies
\begin{alignat*}{1}
&\ \ \ \ \int_{S\cap\{|\nabla u_0| > |\nabla v|\}}\left|\nabla v\right|^2\int^1_0\int^t_0\frac{1}{\left(\left|\nabla u_0\right| + s\left|\nabla v\right|\right)^{2-p}}\,dsdtdx \\
&\geq \int_{S\cap\{|\nabla u_0| > |\nabla v|\}}\frac{\left|\nabla v\right|^2}{|\nabla u_0|^{2-p}}\int^1_0\int^t_0\frac{ds\,dt}{(1+s)^{2-p}}\,dx \\
&= C(p) \int_{S\cap\{|\nabla u_0| > |\nabla v|\}}\frac{\left|\nabla v\right|^2}{|\nabla u_0|^{2-p}}\,dx.
\end{alignat*}
The H\"{o}lder's inequality applied with exponents $\frac{2}{p}$ and $\frac{2}{2-p}$ implies that
\begin{equation*}
\int_{S\cap\{|\nabla u_0| > |\nabla v|\}}\left|\nabla v\right|^p \leq \left(\int_{S\cap\{|\nabla u_0| > |\nabla v|\}}\frac{|\nabla v|^2}{|\nabla u_0|^{2-p}}\right)^{\frac{p}{2}}\left(\int_{S\cap\{|\nabla u_0| > |\nabla v|\}}|\nabla u_0|^p\right)^{\frac{2-p}{2}},
\end{equation*}
or equivalently
\begin{alignat*}{1}
\int_{S\cap\{|\nabla u_0| > |\nabla v|\}}\frac{|\nabla v|^2}{|\nabla u_0|^{2-p}} &\geq \frac{\left(\int_{S\cap\{|\nabla u_0| > |\nabla v|\}}\left|\nabla v\right|^p\right)^{\frac{2}{p}}}{\left(\int_{S\cap\{|\nabla u_0| > |\nabla v|\}}|\nabla u_0|^p\right)^{\frac{2-p}{p}}} \\
&\geq \frac{\left(\int_{S\cap\{|\nabla u_0| > |\nabla v|\}}\left|\nabla v\right|^p\right)^{\frac{2}{p}}}{\left(\int_{\Omega}|\nabla u_0|^p\right)^{\frac{2-p}{p}}}.
\end{alignat*}
Consequently,
\begin{alignat*}{1}
I[v] &\geq C(p)\int_{S\cap\{|\nabla u_0|\leq |\nabla v|\}}|\nabla v|^p + C(p)\frac{\left(\int_{S\cap\{|\nabla u_0| > |\nabla v|\}}\left|\nabla v\right|^p\right)^{\frac{2}{p}}}{\left(\int_{\Omega}|\nabla u_0|^p\right)^{\frac{2-p}{p}}} \\
&\geq C(p)\left(\int_{S\cap\{|\nabla u_0|\leq |\nabla v|\}}|\nabla v|^p\right)^{\frac{2}{p}} + C(p)\frac{\left(\int_{S\cap\{|\nabla u_0| > |\nabla v|\}}\left|\nabla v\right|^p\right)^{\frac{2}{p}}}{\left(\int_{\Omega}|\nabla u_0|^p\right)^{\frac{2-p}{p}}},\ \ \text{as $\delta$ is small} \\
&\geq C(p)A(u_0)\left(\left(\int_{S\cap\{|\nabla u_0|\leq |\nabla v|\}}|\nabla v|^p\right)^{\frac{2}{p}} + \left(\int_{S\cap\{|\nabla u_0|> |\nabla v|\}}|\nabla v|^p\right)^{\frac{2}{p}}\right) \\
&\geq C(p)A(u_0)\left(\int_{S}|\nabla v|^p\right)^{\frac{2}{p}} = C(p)A(u_0)\delta^2,
\end{alignat*}
where the last inequality is a consequence of the elementary inequality
\begin{equation*}
a^{\frac{2}{p}} + b^{\frac{2}{p}} \geq C(p)\left(a+b\right)^{\frac{2}{p}}\ \ \text{for\ }\ a, b\geq 0,
\end{equation*}
and the constant
\begin{equation*}
A(u_0) = \min\left\{1,\frac{1}{\left(\int_{\Omega}|\nabla u_0|^p\right)^{\frac{2-p}{p}}}\right\}.
\end{equation*}

So we have proved $I[v] \geq a > 0$ for some $a> 0$ whenever $v\in C^{\infty}_0(\Omega)$ satisfies $\|v\|_{\mathfrak{B}} = \delta$, for any $p\in (1, \infty)$.
\end{pf}

Let
\begin{equation*}
\mathcal{G} = \{\gamma\in C([0,1],H): \gamma(0) = 0\ \text{and\ }\gamma(1) = v_2\}
\end{equation*}
and
\begin{equation*}
c = \inf_{\gamma\in\mathcal{G}}\max_{0\leq t\leq 1}I[\gamma(t)].
\end{equation*}
The verified Palais-Smale condition and the preceding lemma allow us to apply the Mountain Pass Theorem as stated, for example, in \cite{J} to conclude that there is a $v_1\in \mathfrak{B}$ such that $I[v_1] = c$,
and $I'[v_1] = 0$ in $\mathfrak{B}^*$.
That is
\begin{equation*}
\int_{\Omega}\left|\nabla u_1\right|^{p-2}\nabla u_1\cdot\nabla \varphi + Q(x)\beta_{\varepsilon}(u_1)\varphi dx = 0
\end{equation*}
for any $\varphi\in \mathfrak{B} = W^{1,p}_0(\Omega)$, where $u_1 = v_1 + u_0$. So $u_1$ is a weak solution of the problem (\ref{eulereq}) and (\ref{bdrycondition}). In essence, the Mountain Pass Theorem is a way to produce a saddle point solution. Therefore, in general, $u_1$ tends to be an unstable solution in contrast to the stable solutions $u_0$ and $u_2$.

In this section, we have proved the following theorem.
\begin{theorem}
If $\varepsilon << \sigma_m$ and $J_p(u_2) < J_p(u_0)$, then there exists a third weak solution $u_1$ of the problem (\ref{eulereq}) and (\ref{bdrycondition}). Moreover, $J_p(u_1) \geq J_p(u_0) + a$, where $a$ is independent of $\varepsilon$.
\end{theorem}

\section{A Comparison Principle for Evolution}\label{comparison}
In this section, we prove a comparison theorem for the following evolution problem.
\begin{equation}\label{evolution2}
\left\{\begin{array}{ll}
w_t - \plap w + \alpha(x,w) = 0 &\ \text{in\ }\Omega\times(0, T)\\
w(x,t) = \sigma(x) &\ \text{on\ }\partial \Omega\times (0, T)\\
w(x,0) = v_0(x) &\ \text{for\ }x\in\bar{\Omega},\end{array}\right.
\end{equation}
where $T>0$ may be finite or infinite, and $\alpha$ is a continuous function satisfying $0 \leq \alpha(x,w) \leq Kw$ and
\begin{equation*}
\left|\alpha(x,r_2) - \alpha(x,r_1)\right| \leq K\left|r_2 - r_1\right|
\end{equation*}
for all $x\in\Omega$, $r_1$ and $r_2\in\mathbb{R}$, and some $K \geq 0$. Let us introduce the notation $H_pw = w_t - \plap w + \alpha(x,w)$. We recall a weak sub-solution $w\in L^2(0,T; W^{1,p}(\Omega))$ satisfies
\begin{equation*}
\left.\int_Vw\varphi\ \right|_{t_1}^{t_2} + \int^{t_2}_{t_1}\int_V-w\varphi_t + |\nabla w|^{p-2}\nabla w\cdot\nabla\varphi + \alpha(x,w)\varphi \leq 0
\end{equation*}
for any region $V\subset\subset\Omega$ and any test function $\varphi\in L^2(0,T; W^{1,p}(\Omega))$ such that $\varphi_t\in L^2(\Omega\times\mathbb{R}_T)$ and $\varphi\geq 0$ in $\Omega\times\mathbb{R}_T$, where $L^2_0(0,T; W^{1,p}(\Omega))$ is the subset of $L^2(0,T; W^{1,p}(\Omega))$ that contains functions which is equal zero on the boundary of $\Omega\times\mathbb{R}_T$, where $\mathbb{R}_T = [0, T]$. For convenience, we let $\mathfrak{T}_+$ denote this set of test functions in the following.

In particular, it holds that
\begin{equation*}
\int^T_0\int_{\Omega}-w\varphi_t + <|\nabla w|^{p-2}\nabla w, \nabla\varphi> + \alpha(x,w)\varphi \leq 0
\end{equation*}
for any test function $\varphi\in L^2_0(0,T; W^{1,p}(\Omega))$ such that $\varphi_t\in L^2(\Omega\times\mathbb{R}_T)$ and $\varphi\geq 0$ in $\Omega\times\mathbb{R}_T$.

The comparison principle for weak sub- and super-solutions is stated as follows.
\begin{theorem}\label{paraboliccomparison}
Suppose $w_1$ and $w_2$ are weak sub- and super-solutions of the evolutionary problem
(\ref{evolution2}) respectively with $w_1\leq w_2$ on the parabolic boundary $(\bar{\Omega}\times\{0\})
\cup(\partial \Omega\times (0, +\infty))$. Then $w_1\leq w_2$ in $\mathcal{D}$.
\end{theorem}
Uniqueness of a weak solution of (\ref{evolution2}) follows from the comparison principle, Theorem \ref{paraboliccomparison}, immediately.
\begin{lemma}
For $T > 0$ small enough, if $H_pw_1 \leq 0 \leq H_pw_2$ in the weak sense in
$\Omega\times \mathbb{R}_T$ and $w_1 < w_2$ on $\partial_p(\Omega \times
\mathbb{R}_T)$, then $w_1\leq w_2$ in $\Omega\times\mathbb{R}_T$.
\end{lemma}
\begin{pf}
For any given small number $\delta > 0$, we define a new function
$\tilde{w}_1$ by $$\tilde{w}_1(x,t) = w_1(x,t) - \frac{\delta}{T-t},$$ where
$x\in\bar{\Omega}$ and $0\leq t < T$. In order to prove $w_1\leq
w_2$ in $\Omega\times\mathbb{R}_T$, it suffices to prove $\tilde{w}_1\leq w_2$ in
$\Omega\times\mathbb{R}_T$ for all small $\delta > 0$. Clearly, $\tilde{w}_1 <
w_2$ on $\partial_p(\Omega\times \mathbb{R}_T)$, and $\lim_{t\rightarrow
T}\tilde{w}_1(x,t) = -\infty$ uniformly on $\Omega$. Moreover, the following holds for any $\varphi\in\mathfrak{T}_+$:
\begin{alignat*}{1}
&\ \ \ \ \ \ \int^T_0\int_{\Omega}-\tilde{w}_1\varphi_t + <|\nabla\tilde{w}_1|^{p-2}\nabla\tilde{w}_1, \nabla\varphi> + \alpha(x,\tilde{w}_1)\varphi \\
&= \int^T_0\int_{\Omega}-w_1\varphi_t + <|\nabla w_1|^{p-2}\nabla w_1, \nabla\varphi> + \frac{\delta}{T-t}\varphi_t + \left(\alpha(x,\tilde{w}_1) - \alpha(x, w_1)\right)\varphi \\
&\leq \int^T_0\int_{\Omega}\frac{\delta}{T-t}\varphi_t + K\frac{\delta}{T-t}\varphi, \ \ \text{as $w_1$ is a weak sub-solution}\\
&= \int^T_0\int_{\Omega}\left(-\frac{\delta}{(T-t)^2} + K\frac{\delta}{T-t}\right)\varphi \\
&\leq \int^T_0\int_{\Omega} -\frac{\delta}{2(T-t)^2}\varphi,\ \ \text{for $T\leq\frac{1}{2K}$ so that $2K\leq \frac{1}{T-t}$} \\
&< 0,
\end{alignat*}
i.\,e.\,
\begin{equation*}
H_p\tilde{w}_1 \leq -\frac{\delta}{2(T-t)^2} \leq -\frac{\delta}{2T^2} < 0\ \ \text{in the weak sense.}
\end{equation*}

That is, if we abuse the notation a little by denoting $\tilde{w}_1$ by $w_1$ in the following for convenience, it holds for any $\varphi\in\mathfrak{T}_+$,
\begin{equation*}
\int^T_0\int_{\Omega}-w_1\varphi_t + <|\nabla w_1|^{p-2}\nabla w_1, \nabla\varphi> + \alpha(x,w_1)\varphi \leq \int^T_0\int_{\Omega}-\frac{\delta}{2T^2}\varphi < 0.
\end{equation*}
Meanwhile, for any $\varphi\in \mathfrak{T}_+$, $w_2$ satisfies
\begin{equation*}
\int^T_0\int_{\Omega}-w_2\varphi_t + <|\nabla w_2|^{p-2}\nabla w_2, \nabla\varphi> + \alpha(x,w_2)\varphi \geq 0.
\end{equation*}

Define, for $j = 1, 2$, $v_j(x,t) = e^{-\lambda t}w_j(x, t)$, where the constant $\lambda > 2K$. Then $w_j(x,t) = e^{\lambda t}v_j(x, t)$, and it is clear
that $w_1\leq w_2$ in $\Omega\times\mathbb{R}_T$ is equivalent to $v_1\leq v_2$ in $\Omega\times\mathbb{R}_T$. In addition, for any $\varphi\in \mathfrak{T}_+$, the following inequalities hold:
\begin{alignat*}{1}
&\int^T_0\int_{\Omega}-e^{\lambda t}v_1\varphi_t + e^{\lambda(p-1)t}<|\nabla v_1|^{p-2}\nabla v_1, \nabla\varphi> + \alpha(x,e^{\lambda t}v_1)\varphi \leq -\int^T_0\int_{\Omega}\frac{\delta}{2T^2}\varphi \\
&\text{and\ }\ \int^T_0\int_{\Omega}-e^{\lambda t}v_2\varphi_t + e^{\lambda(p-1)t}<|\nabla v_2|^{p-2}\nabla v_2, \nabla\varphi> + \alpha(x,e^{\lambda t}v_2)\varphi \geq 0.
\end{alignat*}
Consequently, it holds for any $\varphi\in \mathfrak{T}_+$
\begin{alignat*}{1}
&\int^T_0\int_{\Omega}-e^{\lambda t}(v_1 - v_2)\varphi_t + e^{\lambda(p-1)t}<|\nabla v_1|^{p-2}\nabla v_1 - |\nabla v_2|^{p-2}\nabla v_2, \nabla\varphi> \\
 &\ \ \ \ \ \ \ \ + \left( \alpha(x, e^{\lambda t}v_1) - \alpha(x, e^{\lambda t}v_2)\right)\varphi \leq -\int^T_0\int_{\Omega}\frac{\delta}{2T^2}\varphi.
\end{alignat*}
We take $\varphi = \left(v_1 - v_2\right)^+ = \max\{v_1 - v_2, 0\}$ as the test function, since it vanishes on the boundary of $\Omega\times\mathbb{R}_T$. Then
\begin{alignat*}{1}
&\int^T_0\int_{\{v_1>v_2\}}-e^{\lambda t}(v_1 - v_2)(v_1 - v_2)_t + e^{\lambda(p-1)t}<|\nabla v_1|^{p-2}\nabla v_1 - |\nabla v_2|^{p-2}\nabla v_2, \nabla v_1 - \nabla v_2> \\
&\ \ \ \ \ \ \ \ + \left(\alpha(x,e^{\lambda t}v_1) - \alpha(x,e^{\lambda t}v_2)\right)(v_1 - v_2) \leq -\frac{\delta}{2T^2}\int^T_0\int_{\{v_1>v_2\}}(v_1-v_2).
\end{alignat*}
Since
\begin{equation*}
\{v_1>v_2\}\subset\Omega\times(0,T)\ \ \text{due to the facts $v_1\leq v_2$ on $\partial_p(\Omega\times\mathbb{R}_T)$ and $v_1\rightarrow -\infty$ as $t\uparrow T$,}
\end{equation*}
the divergence theorem implies
\begin{equation*}
\int^T_0\int_{\{v_1 > v_2\}}-e^{\lambda t}(v_1-v_2)(v_1-v_2)_t = \int^T_0\int_{\{v_1>v_2\}}\lambda e^{\lambda t}\frac{1}{2}(v_1-v_2)^2.
\end{equation*}
On the other hand,
\begin{equation*}
\left(\alpha(x,e^{\lambda t}v_1) - \alpha(x,e^{\lambda t}v_2)\right)(v_1 - v_2)\geq -Ke^{\lambda t}(v_1-v_2)^2\ \ \text{on $\{v_1 > v_2\}$.}
\end{equation*}
As a consequence, it holds that
\begin{alignat*}{1}
&\int^T_0\int_{\{v_1>v_2\}}\left(\frac{\lambda}{2} - K\right)e^{\lambda t}(v_1 - v_2)^2 + e^{\lambda(p-1)t}<|\nabla v_1|^{p-2}\nabla v_1 - |\nabla v_2|^{p-2}\nabla v_2, \nabla v_1 - \nabla v_2> \\
& \leq -\frac{\delta}{2T^2}\int^T_0\int_{\{v_1>v_2\}}(v_1-v_2).
\end{alignat*}
We call into play two elementary inequalities (\cite{L}) associated with the $p$-Laplacian:
\begin{equation*}
<|b|^{p-2}b - |a|^{p-2}a, b-a> \geq (p-1)|b-a|^2\left(1 + |b|^2 + |a|^2\right)^{\frac{p-2}{2}}\ \ (1\leq p\leq 2),
\end{equation*}
and
\begin{equation*}
<|b|^{p-2}b - |a|^{p-2}a, b-a> \geq 2^{2-p}|b-a|^p\ \ (p\geq 2)\ \ \text{for any $a$, $b\in\mathbb{R}^n$.}
\end{equation*}
By applying them with $b = \nabla v_1$ and $a = \nabla v_2$ in the preceding inequalities, we obtain
\begin{alignat*}{1}
&\int^T_0\int_{\{v_1>v_2\}}\left(\frac{\lambda}{2} - K\right)e^{\lambda t}(v_1 - v_2)^2 + (p-1)e^{\lambda(p-1)t}\left|\nabla v_1 - \nabla v_2\right|^2\left( 1 + |\nabla v_1|^2 + |\nabla v_2|^2\right)^{\frac{p-2}{2}} \\
& \leq -\frac{\delta}{2T^2}\int^T_0\int_{\{v_1>v_2\}}(v_1-v_2)\ \ \ \ \text{for $1 < p <2$}
\end{alignat*}
and
\begin{alignat*}{1}
&\int^T_0\int_{\{v_1>v_2\}}\left(\frac{\lambda}{2} - K\right)e^{\lambda t}(v_1 - v_2)^2 + 2^{2-p}e^{\lambda(p-1)t}\left|\nabla v_1 - \nabla v_2\right|^p \\
& \leq -\frac{\delta}{2T^2}\int^T_0\int_{\{v_1>v_2\}}(v_1-v_2)\ \ \ \ \text{for $p\geq 2$.}
\end{alignat*}
One can easily see in either case the respective inequality is true only if the measure of the set $\{v_1>v_2\}$ is zero. The proof is complete.
\end{pf}

In the next lemma, we show the strict inequality on the boundary data can be relaxed to a non-strict one.
\begin{lemma}
For $T > 0$ sufficiently small, if $H_pw_1 \leq 0 \leq H_pw_2$ in the weak sense in
$\Omega\times\mathbb{R}_T$ and $w_1\leq w_2$ on $\partial_p(\Omega \times
\mathbb{R}_T)$, then $w_1\leq w_2$ on $\overline{\Omega\times\mathbb{R}_T}$.
\end{lemma}
\begin{pf}
For any $\delta > 0$, take $\tilde{\delta} > 0$ such that $\tilde{\delta} \leq\frac{\delta}{4K}$ and define
\begin{equation*}
\tilde{w}_1(x,t) = w_1(x,t) - \delta t - \tilde{\delta}\ \ (x,t)\in\bar{\Omega}\times\mathbb{R}^n.
\end{equation*}
Then $\tilde{w}_1 < w_1 \leq w_2$ on $\partial_p(\Omega\times\mathbb{R}^n)$, and for any $\varphi\in\mathfrak{T}_+$, the following holds:
\begin{alignat*}{1}
&\ \ \ \ \int^T_0\int_{\Omega}-\tilde{w}_1\varphi_t + <|\nabla\tilde{w}_1|^{p-2}\nabla\tilde{w}_1,\nabla\varphi> + \alpha(x,\tilde{w})\varphi \\
&= \int^T_0\int_{\Omega}-w_1\varphi_t + <|\nabla w_1|^{p-2}\nabla w_1,\nabla\varphi> + \alpha(x,w_1)\varphi\\
&\ \ \ \ \ \ \ \ \ \ - \delta \varphi + \left(\alpha(x,w_1 - \delta t - \tilde{\delta}) - \alpha(x,w_1)\right)\varphi \\
&\leq \int^T_0\int_{\Omega} -\delta\varphi + K\left(\delta t + \tilde{\delta}\right)\varphi \leq \int^T_0\int_{\Omega} -\delta\varphi + K\left(\delta T + \tilde{\delta}\right)\varphi \\
&\leq \int^T_0\int_{\Omega}\left(-\delta + \frac{\delta}{2} + \frac{\delta}{4}\right)\varphi\ \ \text{for $T$ small} \\
&= -\frac{\delta}{4}\int^T_0\int_{\Omega}\varphi.
\end{alignat*}
The preceding lemma implies $\tilde{w}_1\leq w_2$ in $\overline{\Omega\times\mathbb{R}_T}$ for small $T$ and for any small $\delta > 0$, and whence the conclusion of this lemma.
\end{pf}

Now the parabolic comparison theorem (\ref{paraboliccomparison}) follows from
the preceding lemma quite easily as shown by the following argument: Let $T_0 > 0$ be
any small value of $T$ in the preceding lemma so that the conclusion of the preceding
lemma holds. Then $w_1\leq w_2$ on $\overline{\Omega\times (0, T_0)}$. In particular,
$w_1 \leq w_2$ on $\partial_p(\Omega\times (T_0,2T_0))$. The preceding lemma may be
applied again to conclude that $w_1 \leq w_2$ on $\overline{\Omega\times (T_0, 2T_0)}$.
And so on. This recursion allows us to conclude that $w_1 \leq w_2$ on $\overline{\Omega\times\mathbb{R}_T}$.


\section{Convergence of Evolution}

Define $\mathfrak{S}$ to be the set of weak solutions of the stationary problem (\ref{eulereq}) and (\ref{bdrycondition}). The $p$-harmonic function $u_0$ is the maximum element in $\mathfrak{S}$, while $u_2$ denotes the least solution which may be constructed as the infimum of super-solutions. We also use the term \textit{non-minimal solution} with the same definition in \cite{CW}. That is, $u$ a non-minimal solution of the problem (\ref{eulereq}) and (\ref{bdrycondition}) if it is a viscosity solution but not a local minimizer in the sense that for any $\delta > 0$, there exists $v$ in the admissible set of the functional $J_p$ with $v = \sigma$ on $\partial\Omega$ such that
$\|v - u\|_{L^{\infty}} < \delta$, and
$J_p(v) < J_p(u)$.

In this section, we consider the evolutionary problem
\begin{equation}\label{evolution}
\left\{\begin{array}{ll}
w_t - \plap w + Q(x)\beta_{\varepsilon}(w) = 0 &\ \text{in\ }\Omega\times(0,+\infty)\\
w(x,t) = \sigma(x) &\ \text{on\ }\partial \Omega\times (0, +\infty)\\
w(x,0) = v_0(x) &\ \text{for\ }x\in\bar{\Omega},\end{array}\right.
\end{equation}
and will apply the parabolic comparison principle (\ref{paraboliccomparison}) proved in Section \ref{comparison} to prove the following convergence of evolution theorem. One just notes that the parabolic problem (\ref{evolution2}) includes the above problem (\ref{evolution}) as a special case so that the comparison principle (\ref{paraboliccomparison}) applies in this case.
\begin{theorem}\label{convergence} If the initial data $v_0$ falls into any of the
categories specified below, the corresponding conclusion of convergence holds.
\begin{enumerate}
\item If $v_0 \leq u_2$ on $\bar{\Omega}$, then $\lim_{t\rightarrow +\infty}w(x,t) = u_2(x)$ locally uniformly
for $x\in\bar{\Omega}$;
\item Define
\begin{equation*}
\bar{u}_2(x) = \inf_{u\in\mathfrak{S}, u \geq u_2, u \neq u_2}u(x),\ x\in\bar{\Omega}.
\end{equation*}
If $\bar{u}_2 > u_2$, then for $v_0$ such that $u_2 < v_0 < \bar{u}_2$, $\lim_{t\rightarrow +\infty}w(x,t) = u_2(x)$ locally uniformly for $x\in\bar{\Omega}$;
\item Define $\bar{u}_0(x) = \sup_{u\in\mathfrak{S}, u \leq u_0, u\neq u_0}u(x)$, $x\in\bar{\Omega}$. If $\bar{u}_0
< u_0$, then for $v_0$ such that $\bar{u}_0 < v_0 < u_0$, $\lim_{t\rightarrow +\infty}w(x,t) =
u_0(x)$ locally uniformly for $x\in\bar{\Omega}$;
\item If $v_0 \geq u_0$ in $\bar{\Omega}$, then $\lim_{t\rightarrow +\infty}w(x,t) = u_0(x)$ locally uniformly
for $x\in\bar{\Omega}$;
\item Suppose $u_1$ is a non-minimal solution of (\ref{eulereq}) and (\ref{bdrycondition}). For any small $\delta > 0$, there exists $v_0$ such that $\|v_0 - u_1\|_{L^{\infty}(\Omega)} < \delta$
and the solution $w$ of the problem (\ref{evolution}) does not satisfy $$\lim_{t\rightarrow
\infty} w(x,t) = u_1(x)\ \text{\ in\ } \Omega.$$
\end{enumerate}
\end{theorem}
\begin{pf}
We first take care of case 4. We may take new initial data a smooth function $\tilde{v}_0$ so that $D^2\tilde{v}_0 < -KI$ and $|\nabla \tilde{v}_0| \geq \delta > 0$ on $\bar{\Omega}$. According to the parabolic comparison principle (\ref{paraboliccomparison}), it suffices to prove the solution $\tilde{w}$ generated by the initial data $\tilde{v}_0$ converges locally uniformly to $u_0$ if we also take $\tilde{v}_0$ large than $v_0$, which can easily be done. So we use $v_0$ and $w$ for the new functions $\tilde{v}_0$ and $\tilde{w}$ without any confusion.

For any $V\subset\subset\Omega$ and any nonnegative function $\varphi$ which is independent of the time variable $t$ and supported in $V$, it holds that
\begin{equation*}
\begin{split}
\int_V|\nabla v_0|^{p-2}\nabla v_0\cdot\nabla\varphi &= \int_V - div\left(|\nabla v_0|^{p-2}\nabla v_0\right)\varphi \\
&\geq \int_V M\varphi\ \ \ \ \text{for some $M = M(n,p,K,\delta) > 0$.}
\end{split}
\end{equation*}
The H\"{o}lder continuity of $\nabla w$ up to $t = 0$ as stated in \cite{DiB}, then implies
\begin{equation*}
\int_V|\nabla w|^{p-2}\nabla w\cdot\nabla\varphi \geq \frac{M}{2}\int_V\varphi
\end{equation*}
for any small $t$ in $(0,t_0)$, and any nonnegative function $\varphi$ which is independent of $t$, supported in $V$ and subject to the condition
\begin{equation}\label{conda}
\frac{\int_V|\nabla \varphi|}{\int_V\varphi}\leq A
\end{equation}
for a fixed constant $A > 0$ and some $t_0 > 0$ dependent on $A$. Then the sub-solution condition on $w$
\begin{equation*}
\left.\int_Vw\varphi\right|_{t=t_2} - \left.\int_Vw\varphi\right|_{t=t_1} + \int^{t_2}_{t_1}\int_V|\nabla w|^{p-2}\nabla w\cdot\nabla\varphi \leq 0
\end{equation*}
implies that
\begin{equation*}
\left.\int_Vw\varphi\right|_{t=t_2} - \left.\int_Vw\varphi\right|_{t=t_1} \leq -\frac{M}{2}(t_2-t_1)\int_V\varphi
\end{equation*}
for any small $t_2 > t_1$ in $(0, t_0)$, and any nonnegative function $\varphi$ which is independent of $t$, supported in $V$ and subject to (\ref{conda}).
In particular, $\left.\int_Vw\varphi\right|^{t_2}_{t_1} \leq 0$ for any nonnegative function $\varphi$ independent of $t$, supported in $V$ and subject to (\ref{conda}). So $$w(x, t_2) \leq w(x, t_1)$$ for any $x\in\Omega$ and $0\leq t_1\leq t_2$. Then the parabolic comparison principle readily implies $w$ is decreasing in $t$ for $t$ in $[0, \infty)$. Therefore $w(x,t)\rightarrow u^{\infty}(x)$ locally uniformly as $t\rightarrow\infty$ and hence $u^{\infty}$ is a solution of (\ref{eulereq}) and (\ref{bdrycondition}). Furthermore, the parabolic comparison principle also implies $w(x,t)\geq u_0(x)$ at any time $t > 0$. Consequently, $u^{\infty} = u_0$ as $u_0$ is the greatest solution of (\ref{eulereq}) and (\ref{bdrycondition}).

Next, we briefly explain the proof for case 1. We may take a new smooth initial data $\tilde{v}_0$ such that $\tilde{v}_0$ is very large negative, $D^2\tilde{v}_0 \geq KI$ and $|\nabla \tilde{v}_0| \geq\delta$ on $\bar{\Omega}$ for large constant $K>0$ and constant $\delta > 0$. It suffices to prove the solution $\tilde{w}$ generated by the initial data $\tilde{v}_0$ converges to $u_2$ locally uniformly on $\bar{\Omega}$ as $t\rightarrow\infty$. Following a computation exactly parallel to that in case 4, we can prove $w$ is increasing in $t$ in $[0, \infty)$. So $w$ converges locally uniformly to a solution $u^{\infty}$ of (\ref{eulereq}) and (\ref{bdrycondition}). As $u^{\infty}\leq u_2$ and $u_2$ is the least solution of (\ref{eulereq}) and (\ref{bdrycondition}), we conclude $u^{\infty} = u_2$.

In case 2, we may replace $v_0$ by a strict super-solution of $\plap v - Q\beta_{\varepsilon}(v) = 0$ in $\bar{\Omega}$ between $u_2$ and $\bar{u}_2$, by employing the fact that $u_2$ is the infimum of super-solutions of (\ref{eulereq}) and (\ref{bdrycondition}). Using $v_0$ as the initial data, we obtain a solution $w(x,t)$ of (\ref{evolution}). Then one argues as in case 4 that for any $V\subset\subset\Omega$, there exist constants $A > 0$ and $t_0 > 0$ such that for $t_1 < t_2$ with $t_1$, $t_2\in [0, t_0)$, $\int_Vw\varphi\,|^{t_2}_{t_1} \leq 0$ for any nonnegative function $\varphi$ independent of $t$, supported in $V$ and subject to the condition $\frac{\int_V|\nabla\varphi|}{\int_V\varphi} \leq A$. As a consequence, $w(x,t_1) \geq w(x,t_2)\ \ (x\in\Omega)$. Then the parabolic comparison principle implies $w$ is decreasing in $t$ over $[0, +\infty)$. Therefore $w(x,t)$ converges locally uniformly to some function $u^{\infty}$ as $t\rightarrow\infty$ which solves (\ref{eulereq}) and (\ref{bdrycondition}). Clearly $u_2(x) \leq w(x,t) \leq \bar{u}_2(x)$ from which $u_2(x) \leq u^{\infty}(x) \leq \bar{u}_2(x)$ follows. As $w$ is decreasing in $t$ and $v_0\neq \bar{u}_2$, $u^{\infty} \neq \bar{u}_2$. Hence $u^{\infty} = u_2$.

The proof of case 3 is parallel to that of case 2 with the switch of sub- and super-solutions. Hence we skip it.

In case 5, we pick $v_0$ with $\|v_0 - u_1\|_{L^{\infty}} < \delta$ and $J_p(v_0) < J_p(u_1)$. Let $w$ be the solution of (\ref{evolution}) with $v_0$ as the initial data. Clearly, we may change the value of $v_0$ slightly if necessary so that it is not a solution of the equation
$$-\nabla\cdot\left(\left(\varepsilon + |\nabla u|^2\right)^{p/2-1}\nabla u\right) + Q(x)\beta(u) = 0$$ for any small $\varepsilon > 0$.

Let $w^{\varepsilon}$ be the smooth solution of the uniformly parabolic boundary-value problem
\begin{equation*}
\left\{\begin{array}{ll}
w_t - \nabla\cdot\left(\left(\varepsilon + |\nabla w|^2\right)^{p/2-1}\nabla w\right) + Q\beta(w) = 0 &\ \ \text{in $\Omega\times (0, +\infty)$}\\
w(x,t) = \sigma(x) &\ \ \text{on $\partial\Omega\times (0, +\infty)$}\\
w(x,0) = v_0(x) &\ \ \text{on $\bar{\Omega}$.}
\end{array}\right.
\end{equation*}
$w^{\varepsilon}$ converges to $w$ in $W^{1,p}(\Omega)$ for every $t\in [0, \infty)$ as $\varepsilon\rightarrow 0$.

We define the functional
$$J_{\varepsilon, p}(u) = \frac{1}{p}\int_{\Omega}\left(\varepsilon + |\nabla u|^2\right)^{p/2} + Q(x)\Gamma(u)\,dx.$$
It is easy to see that
$$\int^t_0\int_{\Omega}\left(w^{\varepsilon}_t\right)^2 - \nabla\cdot\left(\left(\varepsilon + |\nabla w^{\varepsilon}|^2\right)^{p/2-1}\nabla w^{\varepsilon}\right)w^{\varepsilon}_t + Q\beta(w^{\varepsilon})w^{\varepsilon}_t = 0.$$
As $w^{\varepsilon}_t = 0$ on $\partial\Omega\times (0,\infty)$, one gets
$$\int^t_0\int_{\Omega}\left(w^{\varepsilon}_t\right)^2 + \left(\varepsilon + |\nabla w^{\varepsilon}|^2\right)^{p/2-1}\nabla w^{\varepsilon}\cdot \nabla w^{\varepsilon}_t + Q(x)\Gamma(w^{\varepsilon})_t = 0,$$
which implies
$$\int^t_0\int_{\Omega}\left(w^{\varepsilon}_t\right)^2 + \frac{1}{p}\left(\left(\varepsilon + |\nabla w^{\varepsilon}|^2\right)^{p/2}\right)_t + Q(x)\Gamma(w^{\varepsilon})_t = 0.$$
Consequently, it holds
\begin{equation*}
\begin{split}
&\ \ \ \ \int^t_0\int_{\Omega}\left(w^{\varepsilon}_t\right)^2 + \frac{1}{p}\int_{\Omega}\left(\varepsilon + |\nabla w^{\varepsilon}(x,t)|^2\right)^{p/2} + Q\Gamma(w^{\varepsilon}(x,t)) \\
&= \frac{1}{p}\int_{\Omega}\left(\varepsilon + |\nabla w^{\varepsilon}(x,0)|^2\right)^{p/2} + Q\Gamma(w^{\varepsilon}(x,0))
\end{split}
\end{equation*}
i.\,e.\,
\begin{equation*}
\int^t_0\int_{\Omega}\left(w^{\varepsilon}_t\right)^2 + J_{\varepsilon, p}(w^{\varepsilon}(\cdot,t)) = J_{\varepsilon, p}(w^{\varepsilon}(\cdot,0)).
\end{equation*}
Therefore
$$J_{\varepsilon, p}(w^{\varepsilon}(\cdot, t) \leq J_{\varepsilon, p}(v_0),$$
which in turn implies
$$J_p(w(\cdot, t) \leq J_p(v_0) < J_p(u_1).$$
In conclusion, $w$ does not converge to $u_1$ as $t\rightarrow\infty$.
\end{pf}

\end{document}